\documentclass[11pt]{article}
\usepackage{amsmath,latexsym,epsf,mathabx}
\usepackage{changepage}
\usepackage[all]{xy}
\usepackage{hyperref}
\usepackage{color}
\begin{document}

\newcommand{\nc}{\newcommand}
\def\PP#1#2#3{{\mathrm{Pres}}^{#1}_{#2}{#3}\setcounter{equation}{0}}
\def\ns{$n$-star}\setcounter{equation}{0}
\def\nt{$n$-tilting}\setcounter{equation}{0}
\def\Ht#1#2#3{{{\mathrm{Hom}}_{#1}({#2},{#3})}\setcounter{equation}{0}}
\def\qp#1{{${(#1)}$-quasi-projective}\setcounter{equation}{0}}
\def\mr#1{{{\mathrm{#1}}}\setcounter{equation}{0}}
\def\mc#1{{{\mathcal{#1}}}\setcounter{equation}{0}}
\def\HD{\mr{Hom}}
\def\HC{\mr{Hom}_{\mc{C}}}
\def\AdT{\mr{Add}_{\mc{T}}}
\def\adT{\mr{add}_{\mc{T}}}
\def\Kb{\mc{K}^b(\mr{Proj}R)}
\def\kb{\mc{K}^b(\mc{P}_R)}
\def\AdpC{\mr{Adp}_{\mc{C}}}
\def\AdpD{\mr{Adp}_{\mc{D}}}
\newtheorem{Th}{Theorem}[section]
\newtheorem{Def}[Th]{Definition}
\newtheorem{Lem}[Th]{Lemma}
\newtheorem{Pro}[Th]{Proposition}
\newtheorem{Cor}[Th]{Corollary}
\newtheorem{Rem}[Th]{Remark}
\newtheorem{Exm}[Th]{Example}
\newtheorem{Sc}[Th]{}
\def\Pf#1{{\noindent\bf Proof}.\setcounter{equation}{0}}
\def\>#1{{ $\Rightarrow$ }\setcounter{equation}{0}}
\def\<>#1{{ $\Leftrightarrow$ }\setcounter{equation}{0}}
\def\bskip#1{{ \vskip 20pt }\setcounter{equation}{0}}
\def\sskip#1{{ \vskip 5pt }\setcounter{equation}{0}}
\def\mskip#1{{ \vskip 10pt }\setcounter{equation}{0}}
\def\bg#1{\begin{#1}\setcounter{equation}{0}}
\def\ed#1{\end{#1}\setcounter{equation}{0}}
\def\KET{T^{^F\bot}\setcounter{equation}{0}}
\def\KEC{C^{\bot}\setcounter{equation}{0}}

\def\jze{{ \begin{pmatrix} 0 & 0 \\ 1 & 0 \end{pmatrix}}\setcounter{equation}{0}}
\def\hjz#1#2{{ \begin{pmatrix} {#1} & {#2} \end{pmatrix}}\setcounter{equation}{0}}
\def\ljz#1#2{{  \begin{pmatrix} {#1} \\ {#2} \end{pmatrix}}\setcounter{equation}{0}}
\def\jz#1#2#3#4{{  \begin{pmatrix} {#1} & {#2} \\ {#3} & {#4} \end{pmatrix}}\setcounter{equation}{0}}


\title{\bf Recollements and $n$-cotorsion pairs
\footnotetext{
E-mail: weiqingcao@jsnu.edu.cn(Weiqing Cao),~weijiaqun@njnu.edu.cn(Jiaqun Wei), kailywu@163.com(Kaili Wu)}
\footnotetext{
Weiqing cao was supported by the Science Foundation of Jiangsu Normal University (No. 21XFRS024). Jiaqun wei was supported by  the National Natural Science Foundation of China (Grant No. 11771212) and a project funded by the Priority Academic Program Development of Jiangsu Higher
 Education Institutions}}
\smallskip
\author{\small Weiqing Cao$^{a}$, Jiaqun Wei$^{b}$, Kaili Wu$^{c}$\\
  \small    $^{a}$School of Mathematics and Statistics, Jiangsu Normal University,
    Xuzhou 221116,  P.R. China\\
\small $^{b}$ School of Mathematics Sciences,
 Nanjing Normal University, Nanjing \rm210023 P.R. China\\
 \small $^{c}$College of Science, Nanjing Forestry University, Nanjing 210037, P.R. China}

\date{}
\maketitle
\baselineskip 15pt
%
%
\begin{abstract}
\vskip 10pt%
In the present paper, we study the relationships of $n$-cotorsion pairs among three abelian categories in a recollement.
Under certain conditions, we present an explicit construction of gluing of $n$-cotorsion pairs in an abelian category $\mc{D}$ with respect to  $n$-cotorsion pairs in abelian categories $\mc{D}^{'}$, $\mc{D}^{''}$ respectively. On the other hand, we study the construction of $n$-cotorsion pairs in abelian categories $\mc{D}^{'}$, $\mc{D}^{''}$ obtained from $n$-cotorsion pairs in an abelian category $\mc{D}$.

\mskip\





\end{abstract}
%
\vskip 30pt

\section{Introduction}
Cotorsion pairs were invented by Scalce in the category of abelian groups in \cite{Salce}. Subsequently, the concept was readily generalized to any abelian category, triangulated category and even exact category (refer to \cite{Gil,IY,Na,Na2}). Crivei and Torrecillas\cite{CT} presented the $n$-special $\mc{A}$-precovers and $m$-special $\mc{B}$-preenvelopes, and established several conditions under which it is possible to obtain such approximations from an $(m,n)$-cotorsion pair $(\mc{A},\mc{B})$.
Huerta, Mendoza and P\'{e}rez\cite{HMP} present the concept of left and right $n$-cotorsion pairs in abelian category $\mc{D}$. This concept and its dual
generalise the notion of complete cotorsion pairs.  Trlifaj \cite{Trlifaj2} studied infinitely generated tilting and cotilting modules, and characterized the induced tilting and cotilting cotorsion pairs. Actually, cotilting (resp. tilting) cotorsion pairs are left (resp. right) $n$-cotorsion pairs for some $n$.

The recollements of triangulated categories were introduced by Beilinson, Berstein and Deligne \cite{BBD} in order to decompose derived categories of sheaves on topogical spaces with the idea that one triangulated category may be ``glued together" from two others, which plays an important role in representation theory of algebras. Recollements of abelian categories and triangulated categories are
closely related, and they have similar properties in many aspects. Gluing techniques with respect to a recollement of triangulated or abelian categories have been investigated for cotorsion pairs \cite{chen}, torsion pairs \cite{MH}, tilting modules \cite{MZ}, and so on.

It is natural to ask  that whether   an $n$-cotorsion pair can be glued from the other two $n$-cotorsion pair in a recollement of abelian categories. In the present paper, we will give an answer to the question. 
Under certain conditions, we present an explicit construction of gluing of $n$-cotorsion pairs in an abelian category $\mc{D}$ with respect to  $n$-cotorsion pairs in abelian categories $\mc{D}^{'}$, $\mc{D}^{''}$ respectively. On the other hand, we also present the construction of $n$-cotorsion pairs in abelian categories $\mc{D}^{'}$, $\mc{D}^{''}$ obtained from $n$-cotorsion pairs in an abelian category $\mc{D}$.

The paper is organized as follows.
In Section $2$, we recall the definitions of $n$-cotorsion pairs and recollements of abelian categories, and show  some properties of recollements which are useful in our article. In section $3$,
 we glue   $n$-cotorsion pairs in a recollement of abelian categories. In Section $4$, we give some applications of our results.
\vskip 10pt
Let $\mc{D}$ be an abelian category.  Denote by proj$\mc{D}$ and inj$\mc{D}$ the classes of projective and injective objects of $\mc{D}$, respectively.
 Let $\mc{C}\subseteq\mc{D}$ be a class of objects of $\mc{D}$. Given an object $D\in\mc{D}$ and a nonnegative integer $m\geq0$, a $\mc{C}$-resolution of $D$ of length $m$ is an exact sequence $$0\to C_{m}\to \cdots\to C_{1}\to C_{0}\to D\to 0$$
in $\mc{D}$, where $C_{i}\in\mc{C}$ for every interger $0\leq i \leq m$. The resolution dimension of $D$ with respect to $\mc{C}$, denoted resdim$_{\mc{C}}(D)$ is defined as the smallest nonnegative integer $m\geq 0$ such that $D$ has a $\mc{C}$-resolution of length $m$. If such $m$ does not exist, we set resdim$_{\mc{C}}(D):=\infty$. Dually, we have the concepts of $\mc{C}$-coresolutions of $D$ of length $m$ and of coresolution dimension of $D$ with respect to $\mc{C}$, denoted by coresdim$_{\mc{C}}(D)$.

With respect to these two homological dimensions, we shall frequently consider the following classes of objects in $\mc{D}$:
$$\widehat{\mc{C}}_{m}:=\{D\in\mc{D}|~ \text{resdim}_{\mc{C}}(D)\leq m\}$$ and
$$\widecheck{\mc{C}}_{m}:=\{D\in\mc{D}|~ \text{coresdim}_{\mc{C}}(D)\leq m\}.$$

Given two classes of objects $\mc{A}$, $\mc{B}\subseteq\mc{D}$ and an integer $i\geq1$, the notation of Ext$^{i}_{\mc{D}}(\mc{A},\mc{B})=0$ will mean that Ext$^{i}_{\mc{D}}({A},{B})=0$ for any $A\in\mc{A}$ and $B\in\mc{B}$. In the case where $\mc{A}=\{M\}$ or $\mc{B}=\{N\}$, we shall write Ext$^{i}_{\mc{D}}(M,\mc{B})=0$ and Ext$^{i}_{\mc{D}}(\mc{A},N)=0$, respectively. Recall the following notations:
 $$\mc{A}^{\perp_{i}}:=\{N\in\mc{D}|~\mr{Ext}^{i}_{\mc{D}}(\mc{A},N)=0\}$$ and
  $$\mc{A}^{\perp}:=\{N\in\mc{D}|~\mr{Ext}^{i}_{\mc{D}}(\mc{A},N)=0, i>0\}.$$
Dually, we have the notations $^{\perp_{i}}\mc{B}$ and  $^{\perp}\mc{B}$.

\section{Prelimilaries}
Now we recall the notion of $n$-cotorsion pairs.
\begin{Def}\rm{$\mr{\cite[Definition ~2.1]{HMP}}$
Let $\mc{A}$ and $\mc{B}$ be classes of objects in $\mc{D}$. We say that $(\mc{A},\mc{B})$ is a left $n$-cotorsion pair in $\mc{D}$ if the following conditions are satisfied:

(1) $\mc{A}$ is closed under direct summands;

(2) $\mr{Ext}^{i}_{\mc{D}}(\mc{A},\mc{B})=0$ for every $1\leq i\leq n$;

(3) For every object $D\in\mc{D}$, there exists a short exact sequence $0\to K\to A\to D\to 0$ where $A\in\mc{A}$ and $K\in\widehat{\mc{B}}_{n-1}$.}
 \end{Def}

 Dually, define $(\mc{A},\mc{B})$ is a right $n$-cotorsion pair in $\mc{D}$ if $(2)$ above is satisfied, with $\mc{B}$ is closed under direct summands, and if every object of $\mc{D}$ can be embedded into an object of $\mc{B}$ with cokernel in $\widecheck{\mc{A}}_{n-1}$. Finally, $\mc{A}$ and $\mc{B}$ form a $n$-cotorsion pair $(\mc{A},\mc{B})$ if $(\mc{A},\mc{B})$ is both a left $n$-cotorsion pair and right $n$-cotorsion pair in $\mc{D}$.

 \begin{Rem}\rm{
 Note that $1$-cotorsion pairs coincide with the concept of complete cotorsion pairs.}

 \end{Rem}

\begin{Th}$\mathrm{\mr{\cite[Theorem ~2.7]{HMP}}}$
Let $\mc{A}$ and $\mc{B}$ be classes of objects in $\mc{D}$. Then, the following are equivalent:

$(a)$ $(\mc{A},\mc{B})$ is a left $n$-cotorsion pair in $\mc{D}$;

$(b)$ $\mc{A}=\bigcap_{i=1}^{n} {^{\perp_{i}}\mc{B}}$, and for any $D\in\mc{D} $ there is a short exact sequence $0\to K\to A\to D\to 0$ where $A\in\mc{A}$ and $K\in\widehat{\mc{B}}_{n-1}$.

\end{Th}

Dually, we could obtain the equivalent condition of a right $n$-cotorsion pair in $\mc{D}$  and then the equivalent condition of an  $n$-cotorsion pair in $\mc{D}$.

\begin{Th}
Let $\mc{A}$ and $\mc{B}$ be classes of objects in $\mc{D}$. Then the following are equivalent:

$(a)$ $(\mc{A},\mc{B})$ is an $n$-cotorsion pair in $\mc{D}$;

$(b)$ $\mc{A}=\bigcap_{i=1}^{n} {^{\perp_{i}}\mc{B}}$,  $\mc{B}=\bigcap_{i=1}^{n} {\mc{A}^{\perp_{i}}}$, and for any $D\in\mc{D} $ there are two short exact sequences $0\to K\to A\to D\to 0$ and $0\to D\to B\to Z\to 0,$ where $A\in\mc{A}$, $B\in\mc{B}$, $K\in\widehat{\mc{B}}_{n-1}$ and  $Z\in\widecheck{\mc{A}}_{n-1}$.

\end{Th}

\begin{Def}\rm{$\mr{\cite[Definition~4.8]{HMP}}$
A left or right $n$-cotorsion pair $(\mc{A},\mc{B})$ in $\mc{D}$ is hereditary if Ext$^{n+1}_{\mc{D}}(\mc{A},\mc{B})=0$.

}
\end{Def}








We recall the notion of recollement of abelian categories.
\begin{Def}\rm{$\mr{\cite{FP}}$
A recollement, denoted by $(\mc{D}^{'},\mc{D},\mc{D}^{''})$, of the abelian categories $\mc{D}^{'}$, $\mc{D}$ and $\mc{D}^{''}$ is a diagram
\begin{eqnarray*}
\xymatrix{{\mathcal{D^{'}}}\ar[rr]|{i_*}&&\mathcal{D}\ar[rr]|{j^*}\ar[ll]<3ex>|{i^!}\ar[ll]<-3ex>|{i^*}&&\mathcal{D^{''}}\ar[ll]<3ex>|{j_*}\ar[ll]<-3ex>|{j_!}}
\end{eqnarray*}
of abelian categories and additive functors such that

(1) $(i^{\ast},i_{\ast})$, $(i_{\ast},i^{!})$, $(j_{!},j^{\ast})$ and $(j^{\ast},j_{\ast})$ are adjoint pairs,

(2) $i_{\ast}$, $j_{!}$ and $j_{\ast}$ are fully faithful,

(3) $\mr{Im}i_{\ast}=\mr{Ker}j^{\ast}$.}
\end{Def}
We list some properties of recollements of abelian categories (see \cite{FP,P1,pc,pc2}), which will be to be used in the sequel.
\begin{Lem}\label{basic properties}
Let $(\mc{D}^{'},\mc{D},\mc{D}^{''})$ be a recollement of abelian categories.

$(1)$ $i^{\ast}j_{!}=0=i^{!}j_{\ast}$.

$(2)$ The functors $i_{\ast},j^{\ast}$ are exact, $i^{\ast},j_{!}$ are right exact, and $i^{!}, j_{\ast}$ are left exact.

$(3)$ The natural transformations $i^{\ast}i_{\ast}\to 1_{\mc{D}^{'}}$, $1_{\mc{D}^{'}}\to i^{!}i_{\ast}$, $1_{\mc{D}^{''}}\to j^{\ast}j_{!}$,  and $j^{\ast}j_{\ast}\to 1_{\mc{D}^{''}}$ are natural isomorphisms.

$(4)$ Let $D\in\mc{D}$. If $i^{\ast}$ is exact, then we have the exact sequence
$$0\to j_{!}j^{\ast}(D)\to D\to i_{\ast}i^{\ast}(D)\to 0.$$

If $i^{!}$ is exact, then we have the exact sequence
$$0\to i_{\ast}i^{!}(D)\to D\to j_{\ast}j^{\ast}(D)\to 0.$$

$(5)$ If $i^{\ast}$ is exact, then $i^{!}j_{!}=0$; and if $i^{!}$ is exact, then $i^{\ast}j_{\ast}=0$.

$(6)$ If $i^{\ast}$ is exact, then $j_{!}$ is exact.

$(7)$ If $i^{!}$ is exact, then $j_{\ast}$ is exact.
\end{Lem}

\begin{Lem}\label{pre-proj-inj}\rm{$\mr{\cite[Remark ~2.5]{P1}}$}
Let $(\mc{D}^{'},\mc{D},\mc{D}^{''})$ be a recollement of abelian categories.

 $(1)$ If $\mc{D}$ has enough projective objects, then  the functor $i^{\ast }:\mc{D}\to\mc{D}^{'}$ preserves projective objects. In this case, $\mc{D}^{'}$ has enough projective objects and $\mr{proj}\mc{D}^{'}=\mr{add}(i^{\ast}(\mr{proj}\mc{D}))$. Dually, if $\mc{B}$ has enough injective  objects, then the functor $i^{!}:\mc{D}\to \mc{D}^{'}$ preserves injective objects. In this case, $\mc{D}^{'}$ has enough injective objects and $\mr{inj}\mc{D}^{'}=\mr{add}(i^{!}(\mr{inj\mc{D}}))$.

 $(2)$ If $\mc{C}$ has enough projective objects, then the functor $j_{!}:\mc{D}^{'}\to \mc{D}$ preserves projective objects. Dually, if $\mc{C}$ has enough injective objects, then the functor $j_{\ast}:\mc{D}^{'}\to \mc{D}$ preserves injective objects.

 $(3)$ If $j_{\ast}$ is exact and $\mc{D}$ has enough projective objects, then $j^{\ast}$ preserves projective objects. In this case, $\mc{D}^{''}$ has enough projective objects and $\mr{proj}\mc{D}^{''}=\mr{add}(j^{\ast}(\mr{proj}\mc{D}))$.

 $(4)$ If $j_{!}$ is exact and $\mc{D}$ has enough injective objects, then $j^{\ast}$ preserves injective objects. In this case, $\mc{D}^{''}$ has enough injective objects and $\mr{inj}\mc{D}^{''}=\mr{add}(j^{\ast}(\mr{inj}\mc{D}))$.
\end{Lem}

\begin{Pro}\label{EXT}\rm{$\mr{\cite[Proposition ~2.8]{MXZ}}$}
Let $(\mc{D}^{'},\mc{D},\mc{D}^{''})$ be a recollement of abelian categories and $k$ any positive integer.

$(1)$ If $\mc{D}$ has enough projective objects and $i^{\ast}$ is exact, then $$\mr{Ext}^{k}_{\mc{D}^{'}}(i^{\ast}(X), Y)\cong\mr{Ext}^{k}_{\mc{D}}(X,i_{\ast}(Y)).$$

$(2)$ If $\mc{D}^{'}$ has enough projective objects and $i^{!}$ is exact, then $$\mr{Ext}^{k}_{\mc{D}}(i_{\ast}(X), Y)\cong\mr{Ext}^{k}_{\mc{D}^{'}}(X,i^{!}(Y)).$$

$(3)$ If $\mc{D}^{''}$ has enough projective objects and $j_{!}$ is exact, then $$\mr{Ext}^{k}_{\mc{D}}(j_{!}(X), Y)\cong\mr{Ext}^{k}_{\mc{D}^{''}}(X,j^{\ast}(Y)).$$

$(4)$ If $\mc{D}$ has enough projective objects and $j_{\ast}$ is exact, then $$\mr{Ext}^{k}_{\mc{D}^{''}}(j^{\ast}(X), Y)\cong\mr{Ext}^{k}_{\mc{D}}(X,j_{\ast}(Y)).$$

\end{Pro}

\section{ $n$-Cotorsion  pairs in a recollement of abelian categories}

Now, we introduce the recollement of $n$-cotorsion pairs.
\begin{Th}\label{main1}
 Let $(\mc{D}^{'},\mc{D},\mc{D}^{''})$ be a recollement of abelian categories with enough projective objects. Assume that $i^{\ast}$ and $i^{!}$ are exact.
If $(\mc{A}^{'}, B^{'})$ is an $n$-cotorsion pair of $\mc{D}^{'}$, $(\mc{A}^{''}, B^{''})$ is an $n$-cotorsion pair of $\mc{D}^{''}$, denote
$$\mc{A}:=\{D\in\mc{D}| ~j^{\ast}(D)\in \mc{A}^{''} ~and ~ i^{\ast}(D)\in \mc{A}^{'}\},$$
$$\mc{B}:=\{D\in\mc{D}| ~j^{\ast}(D)\in \mc{B}^{''} ~and ~ i^{!}(D)\in \mc{B}^{'}\},$$
then
$(1)$ $\mc{A}^{'}=i^{\ast}(\mc{A})$; $\mc{A}^{''}=j^{\ast}(\mc{A})$; $\mc{B}^{'}=i^{!}(\mc{B})$; $\mc{B}^{''}=j^{\ast}(\mc{B})$;

$(2)$ $$\widecheck{\mc{A}}_{n-1}=\{D\in\mc{D}|~j^{\ast}(D)\in \widecheck{\mc{A}^{''}}_{n-1} ~and ~ i^{\ast}(D)\in \widecheck{\mc{A}^{'}}_{n-1}\},$$
$$\widehat{\mc{B}}_{n-1}=\{D\in\mc{D}|~j^{\ast}(D)\in \widehat{\mc{B}^{''}}_{n-1} ~and ~ i^{!}(D)\in \widehat{\mc{B}^{'}}_{n-1}\};$$

$(3)$ $(\mc{A},\mc{B})$ is an $n$-cotorsion pair in $\mc{D}$. In the case, we call that $(\mc{A},\mc{B})$ is ``glued" by $(\mc{A}^{'}, B^{'})$ and $(\mc{A}^{''}, B^{''})$;

$(4)$ If $(\mc{A}^{'}, B^{'})$ is a hereditary $n$-cotorsion pair of $\mc{D}^{'}$, $(\mc{A}^{''}, B^{''})$ is a hereditary $n$-cotorsion pair of $\mc{D}^{''}$, 
then
 $(\mc{A},\mc{B})$ is a hereditary $n$-cotorsion pair in $\mc{D}$.
 \end{Th}

\Pf. (1) It is obvious that $i^{\ast}(\mc{A})\subseteq\mc{A}^{'}$. For each $A^{'}\in \mc{A}^{'}$, we have $i^{\ast}i_{\ast}(A^{'})\cong A^{'}$ and $j^{\ast}i_{\ast}(A^{'})=0\in\mc{A}^{''}$. So, $i_{\ast}(A^{'})\in \mc{A}$, and then $A^{'}\cong i^{\ast}i_{\ast}(A^{'})\in i^{\ast}(\mc{A})$.
Similarly, $j_{!}(A^{''})\in\mc{A}$ for each $A^{''}\in\mc{A}^{''}$ implies $j^{\ast}(\mc{A})=\mc{A}^{''}$; $i_{\ast}(B^{'})\in\mc{B}$ for each $B^{'}\in\mc{B}^{'}$ induces $\mc{B}^{'}=i^{!}(\mc{B})$; $\mc{B}^{''}=j^{\ast}(\mc{B})$ can be obtained from the fact that $j_{\ast}(B^{''})\in\mc{B}$ for each $B^{''}\in\mc{B}^{''}$.

(2) We need only show  $$\widecheck{\mc{A}}_{n-1}=\{D\in\mc{D}|~j^{\ast}(D)\in \widecheck{\mc{A}^{''}}_{n-1} ~and ~ i^{\ast}(D)\in \widecheck{\mc{A}^{'}}_{n-1}\}.$$
For any $D\in\widecheck{\mc{A}}_{n-1}$, we have an exact sequence
$$ 0\to  D\to A_{0}\to A_{1}\to\cdots\to A_{n-1}\to 0,$$
where $A_{i}\in\mc{A}$, $0\leq i\leq n-1$. Since $j^{\ast}$ is exact, we obtain an exact sequence
$$ 0\to j^{\ast}(D)\to j^{\ast}(A_{0})\to j^{\ast}(A_{1})\to\cdots\to j^{\ast}(A_{n-1})\to 0.$$
It follows that $j^{\ast}(A_{i})\in j^{\ast}(\mc{A})=\mc{A}^{''}$, $0\leq i \leq n-1$. So, $j^{\ast}(D)\in \widecheck{\mc{A}^{''}}_{n-1}$. Since $i^{\ast}$ is exact, we obtain an exact sequence
$$ 0\to i^{\ast}(D)\to i^{\ast}(A_{0})\to i^{\ast}(A_{1})\to\cdots\to i^{\ast}(A_{n-1})\to 0. $$
By (1), we obtain $i^{\ast}(D)\in \widecheck{\mc{A}^{'}}_{n-1}$.

Conversely, let $D\in\mc{D}$ satisfy $j^{\ast}(D)\in \widecheck{\mc{A}^{''}}_{n-1}$ and $i^{\ast}(D)\in \widecheck{\mc{A}^{'}}_{n-1}$. Then we have the following two exact sequences
\begin{equation}0\to j^{\ast}(D)\stackrel{f_{0}}\to A^{''}_{0}\stackrel{f_{1}}\to A^{''}_{1}\to\cdots\to A^{''}_{n-1}\to 0\tag{3.1}\end{equation}
and
\begin{equation}0\to i^{\ast}(D)\stackrel{g_{0}}\to A^{'}_{0}\stackrel{g_{1}}\to A^{'}_{1}\to\cdots\to A_{n-1}^{'}\to 0,\tag{3.2}\end{equation}
where $A_{i}^{''}\in \mc{A}^{''}$ and $A^{'}_{i}\in\mc{A}^{'}$, $0\leq i \leq n-1$. Denote $Z_{i}^{''}=\text{Coker}f_{i-1}$ and $Z_{i}^{'}=\text{Coker}g_{i-1}$, where $1\leq i\leq n-2$. Since $i^{\ast}$ is exact, for any $D\in\mc{D}$, there is an exact sequence
$$0\to j_{!}j^{\ast}(D)\to D\to i_{\ast}i^{\ast}(D)\to 0.$$
Applying $j_{!}$  and $i_{\ast}$ to (3.1) and (3.2) respectively, we obtain the following exact sequences
$$0\to j_{!}j^{\ast}(D)\stackrel{j_{!}(f_{0})}\to j_{!}(A^{''}_{0})\stackrel{j_{!}(f_{1})}\to j_{!}(A^{''}_{1})\to\cdots\to j_{!}(A^{''}_{n-1})\to 0$$
and $$0\to i_{\ast}i^{\ast}(D)\stackrel{i_{\ast}(g_{0})}\to i_{\ast}(A^{'}_{0})\stackrel{i_{\ast}(g_{1})}\to  i_{\ast}(A^{'}_{1})\to\cdots\to i_{\ast}(A_{n-1}^{'})\to 0,$$
where $j_{!}(A^{''}_{i})\in\mc{A}$ and $i_{\ast}(A^{'}_{i})\in\mc{A}$, $0\leq i \leq n-1$. Since $i^{\ast}$ and $i^{!}$ are exact, we obtain Ext$^{1}_{\mc{D}}(i_{\ast}i^{\ast}(D),j_{!}(A^{''}_{0}))\cong\text{Ext}^{1}_{\mc{D}^{'}}(i^{\ast}(D),i^{!}j_{!}(A^{''}_{0}))=0$. So, there is a commutative diagram with exact sequences
$$
 \xymatrix{
 & &0\ar[d]&  0 \ar[d]&  0\ar[d]&& \\
&0\ar[r] &j_{!}j^{\ast}(D)\ar[d]\ar[r]&  D \ar[d]\ar[r]&  i_{\ast}i^{\ast}(D)\ar[d]\ar[r]&0& \\
&0\ar[r] &j_{!}(A_{0}^{''})\ar[d]\ar[r]&j_{!}(A_{0}^{''})\oplus i_{\ast}(A_{0}^{'})\ar[d]\ar[r] &  i_{\ast}(A_{0}^{'})\ar[d]\ar[r]&0&  \\
&0\ar[r] &j_{!}(Z_{1}^{''})\ar[r] \ar[d]             & Z_{1}\ar[r]   \ar[d]                 &  i_{\ast}(Z_{1}^{'})\ar[r]\ar[d]&0 \\
& &0& 0& 0&.& \\
  &&&&\\
}
$$
Repeating the process to $0\to j_{!}(Z_{1}^{''})\to Z_{1}\to  i_{\ast}(Z_{1}^{'})\to 0$, finally we obtain $D\in \widecheck{\mc{A}}_{n-1}$.

(3) We firstly claim that $A\in\mc{A}$ if and only if Ext$^{i}_{\mc{D}}(A,\mc{B})=0$, $1\leq i\leq n$. Since $i^{\ast}$ is exact, we have an exact sequence
$$0\to j_{!}j^{\ast}(A)\to A\to i_{\ast}i^{\ast}(A)\to 0.$$
Thus for each $B\in\mc{B}$, we obtain an exact sequence
$$\text{Ext}^{i}_{\mc{D}}(i_{\ast}i^{\ast}(A),B)\to \text{Ext}^{i}_{\mc{D}}(A,B)\to  \text{Ext}^{i}_{\mc{D}}(j_{!}j^{\ast}(A),B).$$
Since $i^{!}$ and $j_{!}$ are exact, we get $\text{Ext}^{i}_{\mc{D}}(i_{\ast}i^{\ast}(A),B)\cong\text{Ext}^{i}_{\mc{D^{'}}}(i^{\ast}(A),i^{!}(B))=0$ and
 $\text{Ext}^{i}_{\mc{D}}(j_{!}j^{\ast}(A),B)\simeq\text{Ext}^{i}_{\mc{D}^{''}}(j^{\ast}(A),j^{\ast}(B))=0$, $1\leq i \leq n$. Thus $\text{Ext}^{i}_{\mc{D}}(A,B)=0$, $1\leq i \leq n$. On the other hand, if Ext$^{i}_{\mc{D}}(A,\mc{B})=0$, $1\leq i\leq n$, then for each $B^{'}\in \mc{B^{'}}$ and $B^{''}\in \mc{B^{''}}$, we have Ext$^{i}_{\mc{D}^{''}}(j^{\ast}(A),B^{''})\cong \text{Ext}^{i}_{\mc{D}}(A,j_{\ast}(B^{''}))=0$
 and Ext$^{i}_{\mc{D}^{'}}(i^{\ast}(A),B^{'})\cong \text{Ext}^{i}_{\mc{D}}(A,i_{\ast}(B^{''}))=0$
 since $j_{\ast}(B^{''})\in\mc{B}$ and $i_{\ast}(B^{''})\in\mc{B}$. This means $j^{\ast}(A)\in\mc{A}^{''}$ and $i^{\ast}(A)\in\mc{A}^{'}$. Therefore, $A\in\mc{A}$.

 Using a similar argument, we could prove that  $B\in\mc{B}$ if and only if Ext$^{i}_{\mc{D}}(\mc{A},{B})=0$, $1\leq i\leq n$.

 For each $D\in\mc{D}$, there exists an exact sequence $0\to K \to A^{''}\to j^{\ast}(D)\to 0 $ in $\mc{D}^{''}$, where $A^{''}\in\mc{A}^{''}$ and $K\in \widehat{\mc{B}^{''}}_{n-1}$. Since $j_{\ast}$ is exact, we have the exact sequence $0\to j_{\ast}(K) \to j_{\ast}(A^{''})\to j_{\ast}j^{\ast}(D)\to 0 $ in $\mc{D}^{''}$ in $\mc{D}$. We obtain the following commutative diagram with exact rows
 $$
 \xymatrix{
&0\ar[r] &j_{\ast}(K)\ar@{=}[d]\ar[r]&  M \ar[d]\ar[r]& D\ar[d]\ar[r]&0& \\
&0\ar[r] &j_{\ast}(K)\ar[r]&j_{\ast}(A^{''})\ar[r] &  j_{\ast}j^{\ast}(D)\ar[r]&0&  \\
}
$$
where $D\to j_{\ast}j^{\ast}(D)$ is the adjunction morphism and $0\to j_{\ast}(K)\to M\to D\to 0$ is an exact sequence
in $\mc{D}$.  This implies that $j^{\ast}(M)\cong j^{\ast}j_{\ast}(A^{''})\cong A^{''}\in\mc{A}^{''}$. Since $i^{\ast}(M)\in\mc{D}^{'}$ and ($\mc{A}^{'},\mc{B}^{'}$) is an $n$-cotorsion pair of $\mc{D}^{'}$, there exists an exact sequence $0\to L\to  A^{'}\to i^{\ast}(M)\to 0$ satisfying $A^{'}\in\mc{A}^{'}$ and $L\in \widehat{\mc{B}^{'}}_{n-1}$. Then $0\to i_{\ast}(L)\to  i_{\ast}(A^{'})\to i_{\ast}i^{\ast}(M)\to 0$  is an exact sequence of $\mc{D}$, and we  get the commutative diagram with exact rows
$$
 \xymatrix{
&0\ar[r] &i_{\ast}(L)\ar@{=}[d]\ar[r]&  N \ar[d]\ar[r]& M\ar[d]\ar[r]&0& \\
&0\ar[r] &i_{\ast}(L)\ar[r]&i_{\ast}(A^{'})\ar[r] &  i_{\ast}i^{\ast}(M)\ar[r]&0&  \\
}
$$
where $M\to  i_{\ast}i^{\ast}(M)$ is the adjunction morphism and $0\to i_{\ast}(L)\to N\to M\to 0 $ is an exact sequence of $\mc{D}$. Thus $i^{\ast}(N)\cong A^{'}\in\mc{A}^{'}$.  Applying $j^{\ast}$ to $0\to i_{\ast}(L)\to N\to M\to 0 $, we have $j^{\ast}(N)\cong j^{\ast}(M)\cong A^{''}\in\mc{A}^{''}$. So $N\in\mc{A}$. Consider the commutative diagram with exact sequences
$$
 \xymatrix{
 &                        &                 &0\ar[d]  & 0\ar[d]&\\
 &0\ar[r]&i_{\ast}(L)\ar[r]\ar@{=}[d]                 &  W \ar[d]\ar[r]                  &j_{\ast}(K) \ar[d]\ar[r] &0 \\
&0\ar[r] &i_{\ast}(L)\ar[r]                        &N\ar[d]\ar[r] &  M\ar[d]\ar[r]&0 &   \\
 &                        &                 &D\ar@{=}[r]  \ar[d]  & D  \ar[d]&\\
 &                        &                 &0  & 0&.\\
}
$$
We claim that $W\in\widehat{\mc{B}}_{n-1}$. Applying $j^{\ast}$  to $0\to i_{\ast}(L)\to W\to j_{\ast}(K)\to 0$, we obtain $j^{\ast}(W)\cong j^{\ast}j_{\ast}(K)\cong K\in\widehat{\mc{B}^{''}}_{n-1}$. Applying $i^{!}$ to $0\to i_{\ast}(L)\to W\to j_{\ast}(K)\to 0$, we obtain $L\cong i^{!}(W)\in\widehat{\mc{B}^{'}}_{n-1}$. Thus $W\in \widehat{\mc{B}}_{n-1}$. Therefore, we get an exact sequence $0\to W\to N\to D\to 0$, where $N\in\mc{A}$ and $W\in \widehat{\mc{B}}_{n-1}$.

Using a similar argument, we could prove that for each $D\in\mc{D}$, there exists an exact sequence $0\to D\to B\to Z\to 0$, where $B\in\mc{B}$ and $Z\in\widecheck{\mc{A}}_{n-1}$.

(4)  For any $A\in\mc{A}$, we obtain an exact sequence $$0\to j_{!}j^{\ast}(A)\to A\to  i_{\ast}i^{\ast}(A)\to 0.$$
Applying  Hom$_{\mc{D}}(-,B)$, $B\in\mc{B}$, we have the exact sequence
$$\mr{Ext}^{n+1}_{\mc{D}}( i_{\ast}i^{\ast}(A),B)\to \mr{Ext}^{n+1}_{\mc{D}}(A,B)\to \mr{Ext}^{n+1}_{\mc{D}}( j_{!}j^{\ast}(A),B).$$
 By the adjointness of $(i_{\ast},i^{!})$
 and $(j_{!},j^{\ast})$, we get $\text{Ext}^{n+1}_{\mc{D}}(i_{\ast}i^{\ast}(A),B)\cong\text{Ext}^{n+1}_{\mc{D^{'}}}(i^{\ast}(A),i^{!}(B))=0$ and
 $\text{Ext}^{n+1}_{\mc{D}}(j_{!}j^{\ast}(A),B)\simeq\text{Ext}^{n+1}_{\mc{D}^{''}}(j^{\ast}(A),j^{\ast}(B))=0$. Thus $\mr{Ext}^{n+1}_{\mc{D}}(A,B)=0$.
 \ \hfill $\Box$

\vskip 10pt
We also could obtain the case for one-side $n$-cotorsion pairs.

\begin{Cor}\label{cor1}
 Let $(\mc{D}^{'},\mc{D},\mc{D}^{''})$ be a recollement of abelian categories with enough projective objects.
Assume that $i^{\ast}$ and $i^{!}$ are exact.
If $(\mc{A}^{'}, B^{'})$ is a left $n$-cotorsion pair of $\mc{D}^{'}$, $(\mc{A}^{''}, B^{''})$ is a left $n$-cotorsion pair of $\mc{D}^{''}$, define

$$\mc{A}:=\{D\in\mc{D}| ~j^{\ast}(D)\in \mc{A}^{''} ~and ~ i^{\ast}(D)\in \mc{A}^{'}\},$$
$$\mc{B}:=\{D\in\mc{D}| ~j^{\ast}(D)\in \mc{B}^{''} ~and ~ i^{!}(D)\in \mc{B}^{'}\},$$
then
 $(\mc{A},\mc{B})$ is a left $n$-cotorsion pair in $\mc{D}$.

 \end{Cor}

When $n=1$, we could obtain the case for complete cotorsion pairs.
\begin{Cor}\label{cor2}
 Let $(\mc{D}^{'},\mc{D},\mc{D}^{''})$ be a recollement of abelian categories with enough projective objects.
Assume that $i^{\ast}$ and $i^{!}$ are exact.
If $(\mc{A}^{'}, B^{'})$ is a  complete cotorsion pair of $\mc{D}^{'}$, $(\mc{A}^{''}, B^{''})$ is a complete cotorsion pair of $\mc{D}^{''}$, define

$$\mc{A}:=\{D\in\mc{D}| ~j^{\ast}(D)\in \mc{A}^{''} ~and ~ i^{\ast}(D)\in \mc{A}^{'}\},$$
$$\mc{B}:=\{D\in\mc{D}| ~j^{\ast}(D)\in \mc{B}^{''} ~and ~ i^{!}(D)\in \mc{B}^{'}\},$$
then 
 $(\mc{A},\mc{B})$ is a complete cotorsion pair in $\mc{D}$.
 \end{Cor}
Now we show that the converse of Theorem \ref{main1} holds true under certain conditions.
\begin{Th}\label{main2}
 Let $(\mc{D}^{'},\mc{D},\mc{D}^{''})$ be a recollement of abelian categories with enough projective objects.
Assume $i^{\ast}$ and $i^{!}$ are exact.
If $(\mc{A}, \mc{B})$ is an $n$-cotorsion pair of $\mc{D}$, satisfying $j_{\ast}j^{\ast}(\mc{B})\subseteq\mc{B}$, $i_{\ast}i^{\ast}(\mc{A})\subseteq\mc{A}$, $i_{\ast}i^{\ast}(\mc{B})\subseteq\mc{B}$, then

$(1)$ $(i^{\ast}(\mc{A}),i^{!}(\mc{B}))$ is an $n$-cotorsion pair of $\mc{D}^{'}$;

$(2)$ $(j^{\ast}(\mc{A}),j^{\ast}(\mc{B}))$ is an $n$-cotorsion pair of $\mc{D}^{''}$;

 $(3)$ If $(\mc{A}_{1},\mc{B}_{1})$ is ``glued" by $n$-cotorsion pairs $(i^{\ast}(\mc{A}),i^{!}(\mc{B}))$ and  $(j^{\ast}(\mc{A}),j^{\ast}(\mc{B}))$, then $\mc{A}_{1}=\mc{A}$ and $\mc{B}_{1}=\mc{B}$;

$(4)$ If $(\mc{A}, \mc{B})$ is a  hereditary $n$-cotorsion pair, so are the $n$-cotorsion pairs in $(1)$ and $(2)$.
\end{Th}

\Pf. (1) We firstly point out that  $i_{\ast}i^{\ast}(\mc{A})\subseteq\mc{A}$ if and only if $i_{\ast}i^{!}(\mc{B})\subseteq \mc{B}$ and that $j_{\ast}j^{\ast}(\mc{B})\subseteq\mc{B}$ if and only if $j_{!}j^{\ast}(\mc{A})\subseteq\mc{A}$. In fact, since $i^{\ast} $, $j_{\ast}$ and $j_{!}$ are exact, for any $A\in\mc{A}$, $B\in\mc{B}$, we have Ext$^{i}_{\mc{D}}(A,i_{\ast}i^{!}(B))\cong \mr{Ext}^{i}_{\mc{D}^{'}}(i^{\ast}(A),i^{!}(B))\cong\mr{Ext}^{i}_{\mc{D}}(i_{\ast}i^{\ast}(A), B)$, and
 Ext$^{i}_{\mc{D}}(A,j_{\ast}j^{\ast}(B))\cong \mr{Ext}^{i}_{\mc{D}^{''}}(j^{\ast}(A),j^{\ast}(B))\cong\mr{Ext}^{i}_{\mc{D}}(j_{!}j^{\ast}(A), B)$, $1\leq i \leq n$.

 Since $i^{\ast}i_{\ast}\cong id_{\mc{D}^{'}}$ and $i_{\ast}i^{!}(\mc{B})\subseteq \mc{B}$, we get $i^{!}(\mc{B})\subseteq i^{\ast}(\mc{B})$. 
Since $i^{!}i_{\ast}\cong id_{\mc{D}^{'}}$ and $i_{\ast}i^{\ast}(\mc{B})\subseteq \mc{B}$, we also get $i^{\ast}(\mc{B})\subseteq i^{!}(\mc{B})$. Then $i^{!}(\mc{B})= i^{\ast}(\mc{B})$.

Secondly, we show that $A^{'}\in i^{\ast}(\mc{A})$ if and only if Ext$^{i}_{\mc{D}}(A^{'},i^{!}(\mc{B}))=0$, $1\leq i \leq n$. On one hand, if $A^{'}\in i^{\ast}(\mc{A})$, there exists $A\in\mc{A}$ such that $A^{'}=i^{\ast}(A)$. Then for any $i^{!}(B)\in i^{!}(\mc{B})$, we get Ext$^{i}_{\mc{D}^{'}}(A^{'},i^{!}(B))=\mr{Ext}^{i}_{\mc{D}^{'}}(i^{\ast}(A),i^{!}(B))\cong\mr{Ext}^{i}_{\mc{D}}(A,i_{\ast}i^{!}(B))=0$, $1\leq i \leq n$. Therefore, $i_{\ast}i^{!}(\mc{B})\subseteq\mc{B}$ implies  Ext$^{i}_{\mc{D}^{'}}(A^{'},i^{!}(B))=0$, $1\leq i \leq n$. On the other hand, if  Ext$^{i}_{\mc{D}^{'}}(A^{'},i^{!}(\mc{B}))=0$, then for any $i^{!}(B)\in i^{!}(\mc{B})$, we have $0=\mr{Ext}^{i}_{\mc{D}^{'}}(A^{'},i^{!}(B))\cong\mr{Ext}^{i}_{\mc{D}}(i_{\ast}(A^{'}),B)$. Thus $i_{\ast}(A^{'})\in\mc{A}$ and then $A^{'}\cong i^{\ast}i_{\ast}(A^{'})\in i^{\ast}(\mc{A})$.

Similarly, we can prove that $B^{'}\in i^{!}(\mc{B})$ if and only if Ext$^{i}_{\mc{D^{'}}}(i^{\ast}(\mc{A}), B^{'}))=0$, $1\leq i \leq n$.

Finally, for any $D^{'}\in\mc{D^{'}}$, there is an exact sequence $0\to K\to A\to i_{\ast}(D^{'})\to 0$, where $A\in\mc{A}$ and $K\in \widehat{\mc{B}}_{n-1}$. Since $i^{\ast }$ is  exact,  $0\to i^{\ast}(K)\to i^{\ast}(A)\to D^{'}\to 0$ is an exact sequence of $\mc{D}^{'}$. Since $K\in \widehat{\mc{B}}_{n-1}$, there is an exact sequence $0\to B_{n-1}\to \cdots\to B_{1}\to B_{0}\to K\to 0$, where $B_{i}\in\mc{B}$, $0\leq i \leq n-1$. Applying $i^{\ast}$ to this sequence, we obtain the exact sequence $0\to i^{\ast}( B_{n-1})\to \cdots\to i^{\ast}(B_{1})\to i^{\ast}(B_{0})\to i^{\ast}(K)\to 0$, where $i^{\ast}(B_{i})\in i^{\ast}(\mc{}B)= i^{!}(\mc{B})$. Thus the exact sequence $0\to i^{\ast}(K)\to i^{\ast}(A)\to D^{'}\to 0$ satisfies $i^{\ast}(A)\in i^{\ast}({\mc{A}})$ and $i^{\ast}(K)\in \widehat{i^{!}(\mc{B})}_{n-1}$. Similarly, we can prove that,  for any $D^{'}\in\mc{D^{'}}$, there is an exact sequence $0\to D^{'}\to V\to Z\to 0$, where $V\in i^{!}(\mc{B})$ and $Z\in\widecheck{i^{\ast}(\mc{A}})_{n-1}$.

(2) Using the similar arguments as (1), we also can obtain $(j^{\ast}(\mc{A}),j^{\ast}(\mc{B}))$ is an $n$-cotorsion pair of $\mc{D}^{''}$.

(3) By definition, $$\mc{A}_{1}=\{D\in\mc{D}|~j^{\ast}(D)\in j^{\ast}(\mc{A})~\text{and}~i^{\ast}(D)\in i^{\ast}(\mc{A})\};$$
$$\mc{B}_{1}=\{D\in\mc{D}|~j^{\ast}(D)\in j^{\ast}(\mc{B})~\text{and}~i^{!}(D)\in i^{!}(\mc{B})\}.$$
It is obvious that $\mc{A}\subseteq\mc{A}_{1}$ and $\mc{B}\subseteq\mc{B}_{1}$. For any $A_{1}\in\mc{A}_{1}$, there is an exact sequence
$0\to j_{!}j^{\ast}(A_{1})\to A_{1}\to i_{\ast}i^{\ast}(A_{1})\to 0$. Thus for any $B\in\mc{B}$, we get an exact sequence $\mr{Ext}^{i}_{\mc{D}}(i_{\ast}i^{\ast}(A_{1}),B)\to \mr{Ext}^{i}_{\mc{D}}(A_{1},B) \to \mr{Ext}^{i}_{\mc{D}}(j_{!}j^{\ast}(A_{1}),B)$. This implies Ext$^{i}_{\mc{D}}(A_{1},B)=0$, since $\mr{Ext}^{i}_{\mc{D}}(i_{\ast}i^{\ast}(A_{1}),B)\cong\mr{Ext}^{i}_{\mc{D}^{'}}(i^{\ast}(A_{1}),i^{!}(B))=0$ and
$\mr{Ext}^{i}_{\mc{D}}(j_{!}j^{\ast}(A_{1}),B)\cong\mr{Ext}^{i}_{\mc{D}^{''}}(j^{\ast}(A_{1}),j^{\ast}(B))=0$, $1\leq i \leq n$,  $1\leq i \leq n$. Therefore, $\mc{A}_{1}\subseteq\mc{A}$. Similarly, we can prove $\mc{B}_{1}=\mc{B}$.

(4) Since $i_{\ast}i^{\ast}(A)\subseteq\mc{A}$  and $j_{\ast}j^{\ast}(B)\subseteq \mc{B}$, $\mr{Ext}^{n+1}_{\mc{D}^{'}}(i^{\ast}(A),i^{!}(B))\cong\mr{Ext}^{n+1}_{\mc{D}}(i_{\ast}i^{\ast}(A), B)=0$ and Ext$^{n+1}_{\mc{D}}(A,j_{\ast}j^{\ast}(B))\cong \mr{Ext}^{n+1}_{\mc{D}^{''}}(j^{\ast}(A),j^{\ast}(B))$ .
 \ \hfill $\Box$

\vskip 10pt

The following are the versions for the recollement of one-side $n$-cotorsion pairs.

\begin{Cor}\label{left}
 Let $(\mc{D}^{'},\mc{D},\mc{D}^{''})$ be a recollement of abelian categories with enough projective objects.
Assume that $i^{\ast}$ and $i^{!}$ are exact.
 Let $(\mc{A}, \mc{B})$ be a left $n$-cotorsion pair of $\mc{D}$.

 $(1)$ If $i_{\ast}i^{!}(\mc{B})\subseteq\mc{B}$ and $i_{\ast}i^{\ast}(\mc{B})\subseteq\mc{B}$ then
  $(i^{\ast}(\mc{A}),i^{!}(\mc{B}))$ is a left $n$-cotorsion pair of $\mc{D}^{'}$;

$(2)$ If  $j_{!}j^{\ast}(\mc{A})\subseteq\mc{A}$ or $j_{\ast}j^{\ast}(\mc{B})\subseteq\mc{B}$, then $(j^{\ast}(\mc{A}),j^{\ast}(\mc{B}))$ is a left $n$-cotorsion pair of $\mc{D}^{''}$.

\end{Cor}

\begin{Cor}
 Let $(\mc{D}^{'},\mc{D},\mc{D}^{''})$ be a recollement of abelian categories with enough projective objects.
Assume that $i^{\ast}$ and $i^{!}$ are exact.
 Let $(\mc{A}, \mc{B})$ be a right $n$-cotorsion pair of $\mc{D}$.

 $(1)$ If $i_{\ast}i^{\ast}(\mc{A})\subseteq\mc{A}$, $i_{\ast}i^{\ast}(\mc{B})\subseteq\mc{B}$ then
  $(i^{\ast}(\mc{A}),i^{!}(\mc{B}))$ is a right $n$-cotorsion pair of $\mc{D}^{'}$;

$(2)$ If  $j_{!}j^{\ast}(\mc{A})\subseteq\mc{A}$ or $j_{\ast}j^{\ast}(\mc{B})\subseteq\mc{B}$, then $(j^{\ast}{(\mc{A})},j^{\ast}(\mc{B}))$ is a right $n$-cotorsion pair of $\mc{D}^{''}$.

\end{Cor}

When $n=1$, we could obtain the case for complete cotorsion pairs.
\begin{Cor}\label{cor4}
 Let $(\mc{D}^{'},\mc{D},\mc{D}^{''})$ be a recollement of abelian categories with enough projective objects.
Assume that $i^{\ast}$ and $i^{!}$ are exact.
If $(\mc{A}, \mc{B})$ is a complete cotorsion pair of $\mc{D}$, satisfying $j_{\ast}j^{\ast}(\mc{B})\subseteq\mc{B}$, $i_{\ast}i^{\ast}(\mc{A})\subseteq\mc{A}$, $i_{\ast}i^{\ast}(\mc{B})\subseteq\mc{B}$, then

$(1)$ $(i^{\ast}(\mc{A}),i^{!}(\mc{B}))$ is a complete cotorsion pair of $\mc{D}^{'}$;

$(2)$ $(j^{\ast}(\mc{A}),j^{\ast}(\mc{B}))$ is a complete cotorsion pair of $\mc{D}^{''}$;

 $(3)$ If $(\mc{A}_{1},\mc{B}_{1})$ is ``glued" by complete cotorsion pairs $(i^{\ast}(\mc{A}),i^{!}(\mc{B}))$ and  $(j^{\ast}(\mc{A}),j^{\ast}(\mc{B}))$, then $\mc{A}_{1}=\mc{A}$ and $\mc{B}_{1}=\mc{B}$.
\end{Cor}

\section{Applications and examples}

~~~~~In this section, we apply the our main results to $n$-tilting modules and $(m,n)$-cotorsion pairs.

Let $R$ be a ring and  mod$R$  be  the category of finitely generated right $R$-modules.
We recall the definition of   $n$-tilting modules.
\begin{Def}\label{cotilting}\rm{\cite{Happel}\cite{Miyashita}
Let $T\in\mr{mod}R$.  $T$ is called $n$-tilting provided that

$(\mr{T}1)$ The projective dimension of $T$ is at most $n$;

$(\mr{T}2)$ $\mr{Ext}^{i}_{R}(T,T)=0$ for any $1\leq i\leq n$;

$(\mr{T}3)$ There  admits an exact sequence
$$0\to R\to T_{0}\to \cdots\to T_{n}\to  0$$
where $T_{i}\in \mr{add}T$.}
\end{Def}

An $n$-cotilting module $C$ is defined dually.
 By \cite{Trlifaj2}, the pairs  $(\widecheck{(\text{add}T)}_{n}, T^{\perp})$ and  $({^{\perp}C}, \widehat{(\text{add}C)}_{n})$  are hereditary and complete cotorsion pairs.

\begin{Lem} 
Let $T$ be an $n$-tilting $R$-module. Then $(\mr{add}T, T^{\perp})$ is a hereditary right $(n+1)$-cotorsion pair.
\end{Lem}



Now, we show the recollments of $n$-tilting modules under certain conditions,
\begin{Th}\label{main3}
 Let $(\mr{mod}R^{'},\mr{mod}{R},\mr{mod}R^{''})$ be a recollement of module categories. Assume that $i^{\ast}$ and $i^{!}$ are exact.
 Let $T^{'}$ and $T^{''}$ be $n$-tilting modules in mod$R^{'}$ and mod$R^{''}$ respectively. Then $j_{!}(T^{''})\oplus i_{\ast}(T^{'})$ is an $n$-tilting $R$-module.
 \end{Th}

\Pf. Let $T=j_{!}(T^{''})\oplus i_{\ast}(T^{'})$. Since $j_{!}$ and $i_{\ast}$ are  exact and preserve projective modules, $T$ satisfies  $(\mathrm{T}1)$. By Proposition \ref{EXT}, it is easy to prove  that  $T$ satisfies  $(\mathrm{T}2)$. 
 By Lemma \ref{basic properties}\ref{pre-proj-inj}, there is an exact sequence $0\to j_{!}j^{\ast}(R)\to R\to i_{\ast}i^{\ast}(R)\to 0$, where $j^{\ast}(R)\in \text{proj}R^{''}$, $i^{\ast}(R)\in \text{proj} R^{'}$. Since $T^{'}$ and $T^{''}$ are $n$-tilting modules in mod$R^{'}$ and mod$R^{''}$ respectively, $j^{\ast}(R)\in\widehat{(\text{add}T^{''})}_{n}$, $i^{\ast}(R)\in \widehat{(\text{add}T^{'})}_{n}$. Then  $j_{!}j^{\ast}(R)\in\widehat{(\text{add}j_{!}(T^{''}))}_{n}$, $i_{\ast}i^{\ast}(R)\in \widehat{(\text{add}i_{\ast}(T^{'}))}_{n}$. By Proposition \ref{EXT}, we could obtain $R\in\widehat{(\text{add}T)}_{n}$.  So, $T$ is an $n$-tilting $R$-module.
 \ \hfill $\Box$

We also could show the recollments of $n$-cotilting modules under certain conditions.
\begin{Th}\label{main4}
 Let $(\mr{mod}R^{'},\mr{mod}{R},\mr{mod}R^{''})$ be a recollement of module categories. Assume that $i^{\ast}$ and $i^{!}$ are exact.
 Let $C^{'}$ and $C^{''}$ be $n$-cotilting modules in mod$R^{'}$ and mod$R^{''}$ respectively. Then $j_{\ast}(C^{''})\oplus i_{\ast}(C^{'})$ is an $n$-cotilting $R$-module.
 \end{Th}

Another generalization of cotorsion pairs is called ($m,n$)-cotorsion pairs.

\begin{Def}\rm{$\mr{\cite[Definition~3.12]{CT}}$
A pair $(\mc{A}, \mc{B})$ of classes of objects in $\mc{D}$ is called an $(m,n)$-cotorsion pair if $\mc{A}^{\perp_{n}}=\mc{B}$ and $^{\perp_{m}}\mc{B}=\mc{A}$.}
\end{Def}

\begin{Th}\label{main4}
 Let $(\mc{D}^{'}, \mc{D},  \mc{D}^{''})$ be a recollement of abelian categories with enough injectives and enough projectives.
Assume $i^{\ast}$ and $i^{!}$ are exact.
If $(\mc{A}, \mc{B})$ is an $(m,n)$-cotorsion pair of $\mc{D}$, satisfying $j_{\ast}j^{\ast}(\mc{B})\subseteq\mc{B}$, $i_{\ast}i^{\ast}(\mc{A})\subseteq\mc{A}$, then

$(1)$ $(i^{\ast}(\mc{A}),i^{!}(\mc{B}))$ is an  $(m,n)$-cotorsion pair of $\mc{D}^{'}$;

$(2)$ $(j^{\ast}(\mc{A}),j^{\ast}(\mc{B}))$ is an  $(m,n)$-cotorsion pair of $\mc{D}^{''}$.

\end{Th}

\Pf.   (1) Take $B^{'}\in {i^{\ast}(\mc{A})}^{\perp_{n}}$. By Proposition \ref{EXT}, we obtain $0=\text{Ext}^{n}_{\mc{D}^{'}}(i^{\ast}(A),B^{'})\simeq\text{Ext}^{n}_{\mc{D}}(A,i_{\ast}(B^{'}))$ for $A\in\mc{A}$. Then $i_{\ast}(B^{'})\in\mc{A}^{\perp_{n}}=\mc{B}$. Then $i^{!}i_{\ast}(B^{'})\simeq B^{'}\in i^{!}(\mc{B})$.
Conversely, take $B^{'}\in i^{!}(\mc{B})$. Then $B^{'}=i^{!}(B)$ for some $B\in\mc{B}$.
Since $i_{\ast}i^{\ast}(\mc{A})\subseteq\mc{A}$, we obtain $\text{Ext}^{n}_{\mc{D}^{'}}(i^{\ast}(A),B^{'})\simeq\text{Ext}^{n}_{\mc{D}^{'}}(i^{\ast}(A),i^{!}(B))\simeq\text{Ext}^{n}_{\mc{D}}(i_{\ast}i^{\ast}(A),B)=0$ for $A\in\mc{A}$. Hence $B^{'}\in i^{\ast}(\mc{A})^{\perp_{n}}$. Then $i^{\ast}(\mc{A})^{\perp_{n}}=i^{!}(\mc{B})$.
Similarly, we could prove ${^{\perp_{m}}i^{!}(\mc{B})}=i^{\ast}(\mc{A})$.

(2) The proof is similar to that of $(1)$.
 \ \hfill $\Box$

\vskip 10pt

We finish with the following example which shows that the condition ``$i^{\ast}$ is exact" in Corollary \ref{left} is not always necessary.
 \begin{Exm}\rm{
Let $R$ be the finite dimensional algebra given by the quiver $1\to 2$ . Define a triangular matrix algebra
 $\Lambda=
 \begin{pmatrix}
 R&0\\
 R&R
 \end{pmatrix}
 $.
 Any module in mod$\Lambda$ can be uniquely written as a triple $(X,Y)_{f}$ with $X,Y\in\mr{mod}R$ and $f\in\mr{Hom}_{R}(X,Y)$. Then $\Lambda$ is a finite dimensional algebra given by the quiver
$$
 
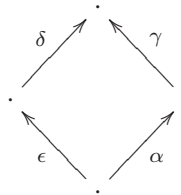
\xymatrix{
 &
                            &  \cdot                                 &  & \\
 &\cdot\ar[ru]^{\delta}    & & \cdot\ar[lu]_{\gamma}& \\
 &                  & \cdot \ar[lu]^{\epsilon} \ar[ru]_{\alpha}                                         & & \\
}
$$
with the relation $\alpha\gamma=\epsilon\delta$. The Auslander-Reiten quiver of $\Lambda$ is
$$
 
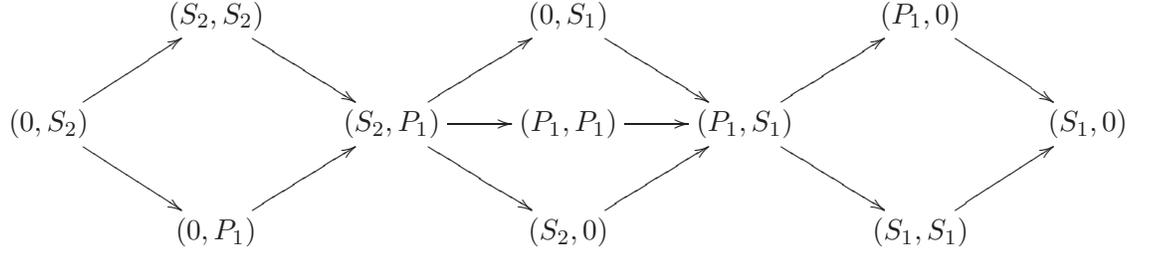
\xymatrix{
 &
                     &                                       & ({S_{2}},{S_{2}})\ar[rd]&&(0,S_{1})\ar[rd] &&({P_{1}},0)\ar[rd]&&\\
 &&({0},{S_{2}})\ar[ru]\ar[rd]& &(S_{2},P_1)\ar[r]\ar[ru]\ar[rd]&  ({P_{1}},{P_{1}})\ar[r]&({P_{1}}, S_1)\ar[ru]\ar[rd]&&({S_{1}},0)&\\
 &&& (0,P_{1}) \ar[ru]&&(S_{2},0)\ar[ru]&&({S_{1}},{S_{1}})\ar[ru]&& \\
}
$$

We have that
\begin{eqnarray*}
\xymatrix{{{\text{mod}R}}\ar[rr]|{i_*}&&\text{mod}\Lambda\ar[rr]|{j^*}\ar[ll]<3ex>|{i^!}\ar[ll]<-3ex>|{i^*}&&\text{mod}{R}\ar[ll]<3ex>|{j_*}\ar[ll]<-3ex>|{j_!}}
\end{eqnarray*}
is a recollement of module categories, where
$$i^{\ast}((X,Y)_f)=\mr{Coker}f, i_{\ast}(Y)=(0,Y), i^{!}((X,Y)_f)=Y;$$
$$j_{!}(Y)=(Y,Y)_1, j^{\ast}((X,Y)_{f})=X, j_{\ast}(X)=(X,0).$$



$(1)$  $(\text{add}R,\mr{mod}R)$ is a $1$-cotorsion pair and $R$ is a $0$-tilting $R$-module. By the construction of Theorem \ref{main3}, $$j_{!}(R)\oplus i_{\ast}(R)=(S(2),S(2))_{1}\oplus (P(1),P(1))_{1}\oplus (0,S(2))\oplus (0,P(1))=\Lambda$$ is a $0$-tilting $\Lambda$-module.
It induces the $1$-cotorsion pair $(\mr{add}\Lambda,\mr{mod}\Lambda)$ in $\mr{mod}\Lambda$. Conversely, by Theorem \ref{main2},  $(\mr{add}\Lambda,\mr{mod}\Lambda)$ induces the $1$-cotorsion pair in mod$R$
$$(i^{\ast}(\mr{add}\Lambda),i^{!}(\mr{mod}\Lambda))=(\mr{add}R,\mr{mod}R)=(j^{\ast}(\mr{add}\Lambda),j^{\ast}(\mr{mod}\Lambda)).$$

$(2)$ It is easy to prove that  $(\mr{add}\Lambda,\mr{add}\Lambda)$ is a left $2$-cotorsion pair in mod$\Lambda$. By Corollary \ref{left}, it induces left $1$-cotorsion pairs in mod$R$
 $$(i^{\ast}(\mr{add}\Lambda),i^{!}(\mr{add}\Lambda))=(\mr{add}R, \mr{add}R)$$ and
$$(j^{\ast}(\mr{add}\Lambda),j^{\ast}(\mr{add}\Lambda))=(\mr{add}R, \mr{add}R).$$}

\end{Exm}

{\small

}

\end{document}